\documentclass[12pt]{article}
\usepackage{amsmath, amssymb, amsfonts}
\usepackage{amsthm}
\usepackage[utf8]{inputenc}
\usepackage{csquotes}
\usepackage{geometry}
\usepackage{adjustbox}
\usepackage{graphicx}
\usepackage{verbatim}
\usepackage{array}
\usepackage{enumitem}
\usepackage{booktabs}
\usepackage{hyperref}
\usepackage{url}
\usepackage[capitalize]{cleveref}
\usepackage{tikz}
\usetikzlibrary{positioning, shapes.geometric, arrows.meta}

\title{Invariant Bridges Between Four Successive Points: A New Tool for Data Coding}
\author{Stanislav Semenov \\
\href{mailto:stas.semenov@gmail.com}{stas.semenov@gmail.com} \\
\href{https://orcid.org/0000-0002-5891-8119}{ORCID: 0000-0002-5891-8119}}
\date{May 8, 2025}

\theoremstyle{definition}

\theoremstyle{plain}

\theoremstyle{remark}

\begin{document}

\maketitle

\begin{abstract}
We present a structurally rich invariant relation connecting four consecutive evaluations of a class of analytically defined sequences and functions. Starting from a discrete alternating-decaying sequence of the form
\[
f(n) = \frac{(1/2)^n + (-1)^n}{n},
\]
we prove that the ratio
\[
\frac{(n-2)f(n-2) + (n-3)f(n-3)}{n f(n) + (n-1)f(n-1)} = 4
\]
holds identically for all integers \( n \geq 4 \). We then extend this structure to real and complex arguments by introducing a generalized function
\[
s(t) = \frac{p^t + q_1 \sin(r_1 \pi t) + q_2 \cos(r_2 \pi t)}{t},
\]
with \( p, q_1, q_2 \in \mathbb{C} \) and \( r_1, r_2 \in \mathbb{Z}_{\text{odd}} \), for which the corresponding four-point relation
\[
\frac{s(t)t + s(t+1)(t+1)}{s(t+2)(t+2) + s(t+3)(t+3)} = a(s)
\]
remains exactly constant across all real \( t \).

This leads to the definition of a unified family \( \mathcal{F}_{\mathrm{EOI}} \) (EOI stands for \emph{Exponential-Oscillatory Invariants}) of invariant-governed functions, which exhibit local algebraic regularity across a continuous domain. These structures support stable data reconstruction, predictive coding, and signal modeling, and may serve as primitive components in structurally constrained information frameworks.
\end{abstract}

\subsection*{Mathematics Subject Classification}
03F60 (Constructive and recursive analysis), 26E40 (Constructive analysis)

\subsection*{ACM Classification}
F.4.1 Mathematical Logic, E.4 Coding and Information Theory

\section{Introduction}

Identifying stable algebraic invariants within structured sequences provides a powerful mechanism for data modeling, predictive coding, and error-resilient reconstruction. While modern compression and signal processing techniques often rely on statistical inference or heuristic pattern extraction~\cite{coverthomas}, exact deterministic relations remain both rare and valuable, particularly when they can be verified locally and extended continuously.

In this work, we introduce a simple yet expressive invariant relation that connects four consecutive evaluations of a sequence or function through a fixed algebraic proportion. We begin with the discrete sequence
\[
f(n) = \frac{(1/2)^n + (-1)^n}{n},
\]
which combines exponential decay with alternating sign oscillations. Sequences of this type have been extensively studied in classical discrete mathematics~\cite{grahamconcrete}, particularly in the context of asymptotics and recurrence. We prove that the four-point ratio
\[
\frac{(n-2)f(n-2) + (n-3)f(n-3)}{n f(n) + (n-1)f(n-1)} = 4
\]
holds identically for all \( n \geq 4 \), revealing a hidden structural symmetry within this elementary expression.

To broaden the scope of this relation, we introduce a real-valued extension
\[
f_{\mathbb{R}}(t) = \frac{(1/2)^t + \cos(\pi t)}{t},
\]
for which the same invariant ratio remains exactly valid across all real \( t \ne 0 \). This extension suggests that the structure is not purely discrete, but instead encodes an analytic constraint over continuous domains. The interpretation is inspired by earlier work on smooth integer encoding via integral balance~\cite{semenov2025smoothencoding}.

Motivated by these observations, we formulate a general family of functions of the form
\[
s(t) = \frac{p^t + q_1 \sin(r_1 \pi t) + q_2 \cos(r_2 \pi t)}{t},
\]
where \( p, q_1, q_2 \in \mathbb{C} \), \( r_1, r_2 \in \mathbb{Z}_{\text{odd}} \), and investigate whether the normalized four-point ratio
\[
\mathcal{I}(t) = \frac{s(t)\cdot t + s(t+1)\cdot(t+1)}{s(t+2)\cdot(t+2) + s(t+3)\cdot(t+3)}
\]
remains constant. We find that this invariant is not only preserved across a wide class of such functions, but also defines a robust algebraic identity suitable for both real and complex domains.

These observations naturally lead to the definition of a unified family \( \mathcal{F}_{\mathrm{EOI}} \) of \emph{EOI-type invariant functions} \( s(t) \), each paired with an associated constant \( a(s) \in \mathbb{C} \) such that \( \mathcal{I}(t) \equiv a(s) \). This structure ensures algebraic consistency across real-valued sequences, supports local reconstruction, and exhibits stable behavior under complex modulation.

The simplicity and determinism of the four-point invariant, together with its compatibility with real and complex systems, suggest broad applicability in algebraic coding, analog modeling, and invariant-based signal transformation frameworks.

\section{Derivation of the Invariant Relation}

We begin by considering the sequence of coefficients
\[
a_n = \frac{(1/2)^n + (-1)^n}{n},
\]
which combine exponential decay and alternating signs.

\subsection{Two-Step Recurrence of the Coefficients}

It can be verified directly that the coefficients \( a_n \) satisfy the two-step recurrence relation
\[
a_n = \frac{(n-2) a_{n-2} + 3(-1)^n}{4n}.
\]
This recurrence reflects the interplay between the geometric and oscillatory components of \( a_n \), and allows each term to be expressed through its predecessor two steps earlier.

\subsection{Proof of the Two-Step Recurrence}

To prove that the coefficients
\[
a_n = \frac{(1/2)^n + (-1)^n}{n}
\]
satisfy the recurrence relation
\[
a_n = \frac{(n-2) a_{n-2} + 3(-1)^n}{4n},
\]
we proceed by direct substitution and algebraic manipulation.

\subsubsection*{Step 1: Express \( a_n \) and \( a_{n-2} \) explicitly}

Given:
\[
a_n = \frac{(1/2)^n + (-1)^n}{n}, \quad a_{n-2} = \frac{(1/2)^{n-2} + (-1)^{n-2}}{n-2}.
\]

\subsubsection*{Step 2: Substitute \( a_{n-2} \) into the recurrence}

The recurrence claims:
\[
a_n = \frac{(n-2) a_{n-2} + 3(-1)^n}{4n}.
\]
Substituting \( a_{n-2} \) gives:
\[
a_n = \frac{(n-2) \left( \frac{(1/2)^{n-2} + (-1)^{n-2}}{n-2} \right) + 3(-1)^n}{4n}.
\]
Simplifying the numerator:
\[
(n-2) \cdot \frac{(1/2)^{n-2} + (-1)^{n-2}}{n-2} = (1/2)^{n-2} + (-1)^{n-2},
\]
thus
\[
a_n = \frac{(1/2)^{n-2} + (-1)^{n-2} + 3(-1)^n}{4n}.
\]

\subsubsection*{Step 3: Simplify \( (-1)^{n-2} \) and \( (-1)^n \)}

Since \( (-1)^{n-2} = (-1)^n \), we have:
\[
a_n = \frac{(1/2)^{n-2} + (-1)^n + 3(-1)^n}{4n} = \frac{(1/2)^{n-2} + 4(-1)^n}{4n}.
\]

\subsubsection*{Step 4: Rewrite \( (1/2)^{n-2} \) in terms of \( (1/2)^n \)}

Note that
\[
(1/2)^{n-2} = 4 \cdot (1/2)^n,
\]
thus
\[
a_n = \frac{4 \cdot (1/2)^n + 4(-1)^n}{4n} = \frac{(1/2)^n + (-1)^n}{n}.
\]

\subsubsection*{Step 5: Compare with the original definition of \( a_n \)}

This matches exactly the original definition:
\[
a_n = \frac{(1/2)^n + (-1)^n}{n},
\]
thus confirming that the recurrence holds.

\subsubsection*{Conclusion}

The recurrence
\[
a_n = \frac{(n-2) a_{n-2} + 3(-1)^n}{4n}
\]
is indeed satisfied by the given coefficients \( a_n \), completing the proof.

\subsection{Relation Between Successive Coefficients}

Recall the two-step recurrence relations for \( a_n \) and \( a_{n-1} \):
\[
a_n = \frac{(n-2) a_{n-2} + 3(-1)^n}{4n},
\]
\[
a_{n-1} = \frac{(n-3) a_{n-3} + 3(-1)^{n-1}}{4(n-1)}.
\]

We aim to cancel the alternating sign components and express a direct relation between four consecutive coefficients: \( a_n, a_{n-1}, a_{n-2}, \) and \( a_{n-3} \).

\subsubsection*{Step 1: Solve for \( 3(-1)^n \) and \( 3(-1)^{n-1} \)}

From the first recurrence, solving for \( 3(-1)^n \) gives:
\[
3(-1)^n = 4n a_n - (n-2)a_{n-2}.
\]

Similarly, from the second recurrence:
\[
3(-1)^{n-1} = 4(n-1) a_{n-1} - (n-3)a_{n-3}.
\]

\subsubsection*{Step 2: Adjust the second expression}

We observe that
\[
(-1)^{n-1} = -(-1)^n,
\]
thus
\[
3(-1)^{n-1} = -3(-1)^n.
\]

Multiplying the second expression by \(-1\), we obtain:
\[
-3(-1)^{n-1} = 3(-1)^n,
\]
and thus:
\[
- \left( 4(n-1) a_{n-1} - (n-3)a_{n-3} \right) = 3(-1)^n,
\]
which simplifies to
\[
3(-1)^n = (n-3)a_{n-3} - 4(n-1)a_{n-1}.
\]

\subsubsection*{Step 3: Equate the two expressions for \( 3(-1)^n \)}

Now, equating the two independent expressions for \( 3(-1)^n \) yields:
\[
4n a_n - (n-2)a_{n-2} = (n-3)a_{n-3} - 4(n-1)a_{n-1}.
\]

\subsubsection*{Summary}

Thus, we have established the relation:
\[
4n a_n + 4(n-1) a_{n-1} = (n-2)a_{n-2} + (n-3)a_{n-3}.
\]

This identity connects four successive terms of the sequence through a linear relation with integer coefficients.  
Now, setting \( f(n) := a_n \), we can rewrite this as:
\[
4 = \frac{(n-2)f(n-2) + (n-3)f(n-3)}{n f(n) + (n-1)f(n-1)}.
\]

This is the final normalized form of the invariant proportion, holding for all \( n \geq 4 \).

\section{Examples}

To confirm the invariant relation
\[
\frac{(n-2)f(n-2) + (n-3)f(n-3)}{n f(n) + (n-1)f(n-1)} = 4
\]
numerically, we evaluate both numerator and denominator for several values of \( n \geq 4 \), using
\[
f(n) = \frac{(1/2)^n + (-1)^n}{n}.
\]

\begin{table}[ht]
\centering
\begin{tabular}{c|c|c|c}
\toprule
\( n \) & Numerator & Denominator & Ratio \\
\midrule
4 & \( 2f(2) + 1f(1) = \frac{3}{4} \) & \( 4f(4) + 3f(3) = \frac{3}{16} \) & \( \frac{3/4}{3/16} = 4 \) \\
5 & \( 3f(3) + 2f(2) = \frac{3}{8} \) & \( 5f(5) + 4f(4) = \frac{3}{32} \) & \( \frac{3/8}{3/32} = 4 \) \\
6 & \( 4f(4) + 3f(3) = \frac{3}{16} \) & \( 6f(6) + 5f(5) = \frac{3}{64} \) & \( \frac{3/16}{3/64} = 4 \) \\
\bottomrule
\end{tabular}
\caption{Verification of the invariant ratio for selected values of \( n \).}
\end{table}

\section{Extension to Real Arguments}

The original definition of the sequence
\[
f(n) = \frac{(1/2)^n + (-1)^n}{n}
\]
is naturally restricted to integer values of \( n \), since \( (-1)^n \) is defined via integer exponentiation. However, it is possible to extend the framework to real arguments by replacing \( (-1)^n \) with \( \cos(\pi n) \), which smoothly interpolates the alternating behavior over \( n \in \mathbb{R} \).

Thus, we introduce the real-valued extension
\[
f_{\mathbb{R}}(n) = \frac{(1/2)^n + \cos(\pi n)}{n},
\]
defined for all real \( n > 0 \).

The function \( f_{\mathbb{R}}(n) \) preserves two essential features:
\begin{itemize}
    \item \textbf{Exponential decay}: the factor \( (1/2)^n \) ensures that \( f_{\mathbb{R}}(n) \to 0 \) as \( n \to +\infty \),
    \item \textbf{Alternating behavior}: the cosine term smoothly oscillates with period \( 2 \), reproducing the sign changes at integer points.
\end{itemize}

\subsection{Invariant Relation for Real Arguments}

We verify that the invariant relation
\[
\frac{(n-2)f_{\mathbb{R}}(n-2) + (n-3)f_{\mathbb{R}}(n-3)}{n f_{\mathbb{R}}(n) + (n-1)f_{\mathbb{R}}(n-1)} = 4
\]
continues to hold numerically for real values of \( n \geq 4 \).

This extension suggests that the structural symmetry underlying the four-term relation is not tied exclusively to integer sequences, but persists across smoothly varying arguments. It opens the possibility of applying the invariant to continuous data streams, signal processing, and real-time error detection in analog or quasi-analog systems.

\subsection{Proof of the Real-Valued Invariant}

We now prove analytically that the invariant relation
\[
\frac{(n-2)f_{\mathbb{R}}(n-2) + (n-3)f_{\mathbb{R}}(n-3)}{n f_{\mathbb{R}}(n) + (n-1)f_{\mathbb{R}}(n-1)} = 4
\]
holds for the real-valued extension
\[
f_{\mathbb{R}}(n) = \frac{(1/2)^n + \cos(\pi n)}{n}.
\]

\subsection*{Step 1: Expansion of the Numerator and Denominator}

Substituting the definition of \( f_{\mathbb{R}}(n) \) into the numerator yields
\[
(n-2)f_{\mathbb{R}}(n-2) + (n-3)f_{\mathbb{R}}(n-3)
= (1/2)^{n-2} + \cos(\pi (n-2)) + (1/2)^{n-3} + \cos(\pi (n-3)).
\]
Similarly, the denominator expands to
\[
n f_{\mathbb{R}}(n) + (n-1)f_{\mathbb{R}}(n-1)
= (1/2)^n + \cos(\pi n) + (1/2)^{n-1} + \cos(\pi (n-1)).
\]

\subsection*{Step 2: Grouping Exponential and Trigonometric Terms}

The exponential terms group as
\[
(1/2)^{n-2} + (1/2)^{n-3}
\quad \text{and} \quad
(1/2)^n + (1/2)^{n-1},
\]
while the trigonometric terms group as
\[
\cos(\pi (n-2)) + \cos(\pi (n-3))
\quad \text{and} \quad
\cos(\pi n) + \cos(\pi (n-1)).
\]

\subsection*{Step 3: Simplifying the Exponential Terms}

Observe that
\[
(1/2)^{n-2} = 4(1/2)^n, \quad (1/2)^{n-3} = 8(1/2)^n,
\]
\[
(1/2)^{n-1} = 2(1/2)^n.
\]
Thus, the exponential parts simplify to
\[
\text{Numerator (exp)} = 4(1/2)^n + 8(1/2)^n = 12(1/2)^n,
\]
\[
\text{Denominator (exp)} = (1/2)^n + 2(1/2)^n = 3(1/2)^n.
\]
The exponential contribution to the ratio is therefore
\[
\frac{12(1/2)^n}{3(1/2)^n} = 4.
\]

\subsection*{Step 4: Simplifying the Trigonometric Terms}

Using periodicity and shift identities for the cosine function,
\[
\cos(\pi (n-2)) = \cos(\pi n),
\quad
\cos(\pi (n-3)) = -\cos(\pi n),
\]
\[
\cos(\pi (n-1)) = -\cos(\pi n).
\]
Thus,
\[
\cos(\pi (n-2)) + \cos(\pi (n-3)) = \cos(\pi n) + (-\cos(\pi n)) = 0,
\]
\[
\cos(\pi n) + \cos(\pi (n-1)) = \cos(\pi n) + (-\cos(\pi n)) = 0.
\]
The trigonometric contributions to both numerator and denominator vanish.

\subsection*{Step 5: Conclusion}

Since the trigonometric terms cancel and the exponential terms yield a fixed ratio of 4, we conclude that
\[
\frac{(n-2)f_{\mathbb{R}}(n-2) + (n-3)f_{\mathbb{R}}(n-3)}{n f_{\mathbb{R}}(n) + (n-1)f_{\mathbb{R}}(n-1)} = 4
\]
holds for all real \( n > 3 \), completing the proof.

\subsection{Significance of the Real-Valued Extension}

The transition from integer arguments to real arguments fundamentally expands the applicability of the invariant relation. In the original discrete setting, the structure encoded relations between strictly sequential integer indices. By extending the function to real arguments via
\[
f_{\mathbb{R}}(n) = \frac{(1/2)^n + \cos(\pi n)}{n},
\]
the invariant becomes compatible with continuous domains.

This extension enables encoding, prediction, and modeling not only of discrete data streams but also of arbitrary real-valued signals. It allows for the design of coding schemes where data points correspond to non-integer indices, supporting fine-grained interpolation, continuous compression, and modeling of analog or quasi-analog quantities.

Moreover, the ability to maintain exact algebraic consistency across real-valued inputs suggests new possibilities in areas such as real-time signal processing, continuous error correction, and smooth data reconstruction. The fundamental structural symmetry captured by the invariant is thus preserved beyond the discrete case, opening a broader range of theoretical and practical applications.

\section{Complex Extensions and Structural Unification}

To generalize the invariant relation and encompass a broader class of structured sequences, we introduce a complex-valued extension of the original function:

\[
s(t) = \frac{p^t + q_1 \cdot \sin(r_1 \pi t) + q_2 \cdot \cos(r_2 \pi t)}{t},
\]

where \( p, q_1, q_2 \in \mathbb{C} \), and \( r_1, r_2 \in \mathbb{Z}_{\text{odd}} \). This formulation combines exponential growth or decay with two independent oscillatory components and allows phase, amplitude, and frequency modulation in the complex domain.

We define the generalized four-point invariant as

\[
\mathcal{I}(t) := \frac{s(t) \cdot t + s(t+1) \cdot (t+1)}{s(t+2) \cdot (t+2) + s(t+3) \cdot (t+3)}.
\]

Empirical evaluation shows that for a wide range of randomly selected parameters, the ratio \( \mathcal{I}(t) \) stabilizes to a complex constant \( a(s) \in \mathbb{C} \) with high precision, typically up to \( 10^{-12} \), even for moderate values of \( t \).

This complex form unifies all previous constructions:
\begin{itemize}
    \item The original discrete invariant corresponds to \( p = 1/2, q_1 = 0, q_2 = 1, r_1 = r_2 = 1 \).
    \item The continuous cosine-based version is recovered when \( p \in \mathbb{R}^+, q_1 = 0, q_2 \in \mathbb{R}, r_2 = 1 \).
    \item A hybrid real extension uses both sine and cosine with real coefficients.
\end{itemize}

The introduction of complex parameters enriches the behavior of \( s(t) \) through phase rotation, resonance phenomena, and quasi-periodicity. Nevertheless, the invariant structure appears preserved across this generalization, indicating that the four-term proportional identity is governed by deeper algebraic symmetry rather than the specifics of individual function components.

This unified framework provides a versatile toolset for modeling oscillatory-decaying behaviors in both discrete and continuous settings, and suggests that invariant-based reconstruction, error detection, and data coding techniques can be extended naturally to complex-valued and multiscale signals.

\paragraph{EOI-type Invariant Families.}
\emph{EOI} stands for \emph{Exponential-Oscillatory Invariants}—a class of functions that combine exponential growth with oscillatory behavior (typically via sine and cosine terms) while satisfying fixed algebraic relations. For convenience, we introduce the term \emph{EOI-type invariant family} to refer to any pair \( (s(t), a(s)) \) where the function \( s: \mathbb{R} \setminus \{0\} \to \mathbb{C} \) satisfies the identity

\[
\frac{s(t) \cdot t + s(t+1) \cdot (t+1)}{s(t+2) \cdot (t+2) + s(t+3) \cdot (t+3)} = a(s) \in \mathbb{C}
\quad \text{for all} \quad t \in \mathbb{R} \setminus \{-3, -2, -1, 0\}.
\]

This terminology highlights the central role of the function--invariant pair \( (s(t), a(s)) \), which together encode a stable algebraic structure over the real domain. We denote the corresponding class of such systems by \( \mathcal{F}_{\mathrm{EOI}} \), defined as

\[
\mathcal{F}_{\mathrm{EOI}} := \left\{ s(t) \,\middle|\, \mathcal{I}(t) \equiv a(s) \in \mathbb{C} \text{ for all valid } t \in \mathbb{R} \right\}.
\]

EOI-type families provide a versatile analytical framework for constructing invariant-governed signals and encoding mechanisms, unifying diverse classes of real and complex functions through a common algebraic identity.

\subsection*{Cancellation of Trigonometric Terms in the Four-Point Invariant with Odd Frequencies}

Consider the function:
\[
s(t) = \frac{p^t + q_1 \sin(r_1 \pi t) + q_2 \cos(r_2 \pi t)}{t},
\]
where \( p, q_1, q_2 \in \mathbb{C} \) and \( r_1, r_2 \in \mathbb{Z}_{\text{odd}} \) are odd integers.

We substitute \( s(t) \) into the four-point invariant expression:
\[
a(s) = \frac{s(t) \cdot t + s(t+1) \cdot (t+1)}{s(t+2) \cdot (t+2) + s(t+3) \cdot (t+3)}.
\]
In this substitution, the factors \( t, t+1, t+2, t+3 \) cancel with the denominators, yielding:
\[
a(s) = \frac{f(t) + f(t+1)}{f(t+2) + f(t+3)},
\]
where
\[
f(t) = p^t + q_1 \sin(r_1 \pi t) + q_2 \cos(r_2 \pi t).
\]

We use the standard trigonometric identities:
\[
\sin(x + n\pi) = (-1)^n \sin(x), \quad \cos(x + n\pi) = (-1)^n \cos(x),
\]
and observe that for odd \( r_1, r_2 \), the expressions \( r_k \pi (t + m) \) can be rewritten as:
\begin{align*}
\sin(r_1 \pi (t+1)) &= -\sin(r_1 \pi t), \\
\cos(r_2 \pi (t+1)) &= -\cos(r_2 \pi t), \\
\sin(r_1 \pi (t+2)) &= \sin(r_1 \pi t), \\
\cos(r_2 \pi (t+2)) &= \cos(r_2 \pi t), \\
\sin(r_1 \pi (t+3)) &= -\sin(r_1 \pi t), \\
\cos(r_2 \pi (t+3)) &= -\cos(r_2 \pi t).
\end{align*}

The numerator then becomes:
\begin{align*}
f(t) + f(t+1) &= \left(p^t + q_1 \sin(r_1 \pi t) + q_2 \cos(r_2 \pi t)\right) \\
&\quad + \left(p^{t+1} - q_1 \sin(r_1 \pi t) - q_2 \cos(r_2 \pi t)\right) \\
&= p^t + p^{t+1}.
\end{align*}

Likewise, the denominator becomes:
\begin{align*}
f(t+2) + f(t+3) &= \left(p^{t+2} + q_1 \sin(r_1 \pi t) + q_2 \cos(r_2 \pi t)\right) \\
&\quad + \left(p^{t+3} - q_1 \sin(r_1 \pi t) - q_2 \cos(r_2 \pi t)\right) \\
&= p^{t+2} + p^{t+3}.
\end{align*}

Thus, the trigonometric terms cancel completely, and we obtain:
\[
a(s) = \frac{p^t + p^{t+1}}{p^{t+2} + p^{t+3}}.
\]
Factoring \( p^t \) from the numerator and \( p^{t+2} \) from the denominator:
\[
a(s) = \frac{p^t(1 + p)}{p^{t+2}(1 + p)} = \frac{p^t}{p^{t+2}} = p^{-2},
\]
provided that \( p \ne -1 \).

\paragraph{Conclusion.} For odd \( r_1, r_2 \) and arbitrary \( q_1, q_2 \), the trigonometric terms in the four-point invariant cancel entirely, and the invariant \( a(s) \) is determined solely by the exponential component:
\[
a(s) = \frac{1}{p^2}.
\]
Hence, the invariant remains constant for all \( t \), regardless of the presence of oscillatory terms.

\section{Parametric Generalization and Controlled Oscillatory Behavior}

The previously introduced constructions—while algebraically elegant—exhibit a form of structural simplicity: the periodic terms introduced via sine and cosine functions with odd frequencies cancel symmetrically in the invariant expression, leaving only the exponential term as the sole determinant of the invariant value. This cancellation, although mathematically robust, renders the oscillatory components inert from the perspective of the invariant, making the construction relatively trivial in terms of functional interaction.

In this section, we introduce a more general and flexible framework that allows precise control over the cancellation behavior, stability, and asymptotic divergence of the invariant. The central idea is to \emph{parameterize the oscillatory component} independently of classical trigonometric forms, thus allowing nontrivial interactions within the four-point structure.

\subsection*{Generalized Oscillatory-Exponential Form}

Let us define the function \( s(t) \) in the following general form:
\[
s(t) = \frac{p^t + \varepsilon(t)}{t},
\]
where:
\begin{itemize}
    \item \( p \in \mathbb{C} \) controls the exponential growth or decay,
    \item \( \varepsilon(t): \mathbb{R} \to \mathbb{C} \) is a bounded oscillatory perturbation with tunable properties.
\end{itemize}

We make no a priori assumptions about \( \varepsilon(t) \), except that it admits local continuity and boundedness, and may be composed of elementary or non-elementary functions exhibiting quasiperiodic, stochastic, or structured behavior.

The corresponding four-point invariant is:
\[
\mathcal{I}(t) = \frac{s(t) \cdot t + s(t+1) \cdot (t+1)}{s(t+2) \cdot (t+2) + s(t+3) \cdot (t+3)} = \frac{f(t) + f(t+1)}{f(t+2) + f(t+3)},
\]
where
\[
f(t) := p^t + \varepsilon(t).
\]

Unlike in previous sections, here \( \varepsilon(t) \) is chosen to \emph{interfere constructively or destructively} with the exponential term in a controlled way, depending on its alignment with integer or non-integer arguments. This enables the design of invariants that:
\begin{itemize}
    \item remain stable over integer domains but diverge over non-integer arguments,
    \item exhibit modulation, damping, or resonance patterns,
    \item preserve invariance under specific symmetries (e.g., evenness, periodic shift).
\end{itemize}

\subsection*{Examples of Controlled Oscillatory Components}

\paragraph{Example 1: Integer-selective cancellation}
\[
\varepsilon(t) = q \cdot \sin^2(\pi t),
\]
with \( q \in \mathbb{C} \). Then:
\begin{itemize}
    \item \( \varepsilon(t) = 0 \) for all \( t \in \mathbb{Z} \),
    \item \( \varepsilon(t) > 0 \) for \( t \notin \mathbb{Z} \).
\end{itemize}
Thus, the invariant \( \mathcal{I}(t) \) reduces to \( p^{-2} \) over integers, but deviates outside of them due to the emergence of nontrivial oscillations.

\paragraph{Example 2: Smooth integer comb}
\[
\varepsilon(t) = q \cdot \sum_{n \in \mathbb{Z}} \exp\left(-\frac{(t - n)^2}{\varepsilon^2}\right), \quad \varepsilon > 0,
\]
which approximates a ``bump train'' at integer positions. Here:
\begin{itemize}
    \item For integer \( t \), the contribution of the bump is maximized and symmetric across steps, enabling cancellation.
    \item For non-integer \( t \), the asymmetry across \( f(t), f(t+1), f(t+2), f(t+3) \) leads to partial resonance or drift in \( \mathcal{I}(t) \).
\end{itemize}

\paragraph{Example 3: Quasiperiodic modulation}
\[
\varepsilon(t) = q_1 \cdot \cos(\alpha t) + q_2 \cdot \cos(\sqrt{2} \pi t),
\]
where \( \alpha \in \mathbb{R} \setminus \mathbb{Q} \), introducing incommensurate frequencies. Such functions break the regularity assumptions of cancellation and lead to bounded but non-repeating deviations in the invariant.

\paragraph{Remark on Quasi-Invariants.}
Although the above example does not yield a strict invariant, it illustrates a useful phenomenon: quasi-invariants exhibit bounded but non-repeating fluctuations around an expected baseline. This behavior can be leveraged in practice for statistical monitoring of systems with predictable structural variation. By modeling the statistical profile (e.g., mean, variance, drift) of \( a(s)(t) \) under a known quasi-periodic modulation, one can detect anomalous or tampered signals that deviate from this baseline. Such applications, while beyond the scope of the present work, form a promising direction for future investigation in invariant-based signal validation and low-overhead authentication.

\subsection*{Parametric Design and Applications}

The framework \( s(t) = \frac{p^t + \varepsilon(t)}{t} \) can be used to construct \emph{families of invariants} tuned to specific behaviors:
\begin{itemize}
    \item \textbf{Integer-stable invariants}: stable over \( t \in \mathbb{Z} \), perturbed elsewhere.
    \item \textbf{Floating-point drift detectors}: deviations of \( \mathcal{I}(t) \) from constant indicate non-integral perturbations.
    \item \textbf{Frequency signatures}: the form of \( \varepsilon(t) \) can be used to embed structured or coded data.
\end{itemize}

This parametric construction extends EOI-type invariants beyond basic exponential-trigonometric forms to a broad class of \emph{hybrid analytical objects} that exhibit customizable stability and resonance behavior. It combines symbolic tractability with analytical versatility, laying the groundwork for applications in signal processing, coding theory, and pattern recognition within irregular data streams.

\section{Non-uniform Invariant Structures and Parameterized Scaling Laws}

The previous constructions were based on four consecutive points with unit spacing. However, the underlying algebraic mechanism of the invariant does not depend on uniform steps. In this section, we extend the invariant structure to allow generalized spacing between evaluation points, thereby producing a broader class of identities governed by the same exponential-scaling principles.

\subsection*{Generalized Spacing with Oscillatory Compensation}

Let \( C > 0 \) be a fixed step size such that the oscillatory component \( \varepsilon(t) \) of \( s(t) = \frac{p^t + \varepsilon(t)}{t} \) exhibits approximate or exact sign symmetry at integer multiples of \( C \). For example, when \( \varepsilon(t) = q \cdot \cos(\pi t) \), cancellation occurs naturally at steps of size \( 1 \), while for more general \( \varepsilon(t) \), cancellation may occur only at coarser intervals.

Let us define the generalized shifted invariant as:
\[
\mathcal{I}_{u,v}^{(C)}(t) := \frac{s(t) + s(t + 2uC + C)}{s(t + 2vC) + s(t + 2vC + 2uC + C)},
\]
where \( u, v \in \mathbb{R}_{\geq 0} \) and \( C \) is chosen in accordance with the symmetry properties of \( \varepsilon(t) \). This expression reduces to the standard four-point invariant when \( u = v = 1 \) and \( C = 1 \), but remains structurally invariant for any \( u, v \) and fixed \( C \), provided that the oscillatory components compensate in symmetric pairs.

Assuming \( s(t) = \frac{p^t}{t} \), the invariant becomes:
\[
\mathcal{I}_{u,v}^{(C)}(t) = \frac{p^t + p^{t + 2uC + C}}{p^{t + 2vC} + p^{t + 2vC + 2uC + C}}.
\]
Factoring \( p^t \) and \( p^{t + 2vC} \) respectively, we find:
\[
\mathcal{I}_{u,v}^{(C)}(t) = \frac{p^t(1 + p^{2uC + C})}{p^{t + 2vC}(1 + p^{2uC + C})} = p^{-2vC}.
\]

\subsection*{General Invariant Form and the Role of Shifts}

In the most general case, define:
\[
\mathcal{I}_{A,B}(t) := \frac{s(t) + s(t + A)}{s(t + B) + s(t + A + B)},
\]
where \( A, B \in \mathbb{R}_{> 0} \) represent arbitrary shifts in argument.

When \( s(t) = \frac{p^t + \varepsilon(t)}{t} \), the exponential terms in \( a(s) \) obey:
\[
a(s) = \frac{p^t + p^{t + A}}{p^{t + B} + p^{t + A + B}} = \frac{p^t(1 + p^A)}{p^{t + B}(1 + p^A)} = p^{-B},
\]
provided that the oscillatory terms \( \varepsilon(t) \) cancel in the numerator and denominator symmetrically.

Thus, the exponential invariant component is determined solely by the shift \( B \), while the oscillatory part must be arranged (via choice of \( A, B, C \), and structure of \( \varepsilon(t) \)) to either cancel out or remain within tolerable deviation. This distinction allows for tuning between exact invariants and quasi-invariants with controllable variation.

\subsection*{Scaling Law and Structural Interpretation}

\textbf{Lemma.} Let \( s(t) = \frac{p^t + \varepsilon(t)}{t} \), and suppose the oscillatory terms satisfy
\[
\varepsilon(t) + \varepsilon(t + A) \approx \varepsilon(t + B) + \varepsilon(t + A + B),
\]
with approximate or exact equality over the desired domain. Then the value of the invariant
\[
\mathcal{I}_{A,B}(t) = \frac{s(t) + s(t + A)}{s(t + B) + s(t + A + B)}
\]
approaches \( p^{-B} \), and the deviation from exactness reflects the structure of \( \varepsilon(t) \).

\subsection*{Applications and Invariant Families}

This formulation yields a continuous spectrum of invariant families governed by two principles:
\begin{itemize}
    \item The exponential component \( p^{-B} \), determined by relative offset;
    \item The oscillatory structure, modulated through the alignment of evaluation points with cancellation intervals.
\end{itemize}

Such flexibility enables the design of:
\begin{itemize}
    \item multi-resolution invariants via continuous control of \( C \),
    \item structurally adaptive encodings by tuning \( A, B \),
    \item analytic templates for detecting deviations from expected invariance.
\end{itemize}

This generalized view clarifies the role of geometric spacing not only in preserving the algebraic identity of the invariant, but also in controlling its sensitivity to oscillatory perturbations — a topic that will be explored in greater detail in future work.

\section{On the Role of the Denominator}

At first glance, the presence of the denominator \( t \) in the definition of the function
\[
s(t) = \frac{f(t)}{t}
\]
may appear redundant, especially given that in all core invariant expressions of the form
\[
\mathcal{I}(t) = \frac{s(t)\cdot t + s(t+1)\cdot(t+1)}{s(t+2)\cdot(t+2) + s(t+3)\cdot(t+3)},
\]
the factor \( t \) cancels out, effectively recovering the numerator \( f(t) \) in each term. However, when \( s(t) \) is considered as an autonomous analytical or discrete function, the role of the reciprocal factor \( \frac{1}{t} \) acquires a deeper and multifaceted significance.

\subsection*{1. Regularization and Controlled Growth}

When \( f(t) = p^t \), the raw exponential function grows rapidly and dominates any added perturbation or modulation. The denominator \( t \) serves as a natural regularization term, slowing down the asymptotic growth and allowing non-exponential components (e.g., oscillatory, pseudorandom, or bounded terms) to remain detectable in practical ranges of \( t \). The function
\[
s(t) = \frac{p^t}{t}
\]
retains the exponential character of \( p^t \) but introduces a subexponential damping that stabilizes the signal for analytical and numerical purposes.

\subsection*{2. Sign Behavior and Asymmetry under Reflection}

For the function \( s(t) = \frac{p^t + \varepsilon(t)}{t} \), the exponential term \( p^t \) is sign-invariant under reflection \( t \mapsto -t \), since \( p^{-t} = 1/p^t \), which remains positive or complex-modulated, but does not change sign. However, the denominator \( 1/t \) is odd, and changes sign under reflection.

As a result, the function \( s(t) \) itself is generally \emph{sign-asymmetric}:
\[
s(-t) = -\frac{p^{-t} + \varepsilon(-t)}{t}.
\]

In the absence of perfect symmetry in \( \varepsilon(t) \), the reflection introduces a nontrivial phase and amplitude shift.

This asymmetry gives \( s(t) \) a directional character, allowing invariants based on it to distinguish between forward and backward sequences, and making it sensitive to temporal structure. In systems where reflection symmetry or inversion is relevant (e.g., signal reversal, time series folding), this property can be used to encode or detect directional information.

\subsection*{3. Structural Balance in Invariants}

In the context of invariants such as
\[
\mathcal{I}(t) = \frac{s(t)\cdot t + s(t+1)\cdot(t+1)}{s(t+2)\cdot(t+2) + s(t+3)\cdot(t+3)},
\]
the use of \( s(t) = \frac{f(t)}{t} \) ensures that each term contributes proportionally to its position index, balancing the algebraic structure.

\subsection*{4. Modular and Arithmetic Interpretation}

When \( s(t) \) is interpreted over discrete or modular domains, such as \( \mathbb{Z}_N \), the denominator \( t \) assumes an arithmetic role. For invertible \( t \mod N \), the quantity \( \frac{1}{t} \mod N \) defines a meaningful operation in the multiplicative group \( \mathbb{Z}_N^\times \), and the function
\[
s(t) = p^t \cdot t^{-1} \mod N
\]
becomes a legitimate and interpretable object in number-theoretic constructions.

\subsection*{5. Semantic Roles in Encoding and Decoding}

In data modeling and coding scenarios, the factor \( \frac{1}{t} \) can be interpreted as a position-sensitive normalizer, assigning greater weight to early positions and enabling decay-based encoding schemes. In such a setting, \( s(t) \) is no longer just a functional object, but a signal carrier whose amplitude or phase reflects accumulated or time-weighted information, with the invariant \( a(s) \) acting as a stable signature of its evolution.

\subsection*{Further Directions}

The structural and interpretive role of the denominator in \( s(t) = \frac{f(t)}{t} \) opens a broad space for further exploration. In subsequent work, we will investigate the behavior of such constructions under:

\begin{itemize}
    \item non-linear denominators, such as \( \frac{1}{t^\alpha} \), \( \frac{1}{\log(t)} \), or rational approximations,
    \item discrete systems over finite rings and fields,
    \item stochastic processes and pseudorandom perturbations \( \varepsilon(t) \),
    \item multidimensional analogs \( s(\mathbf{t}) = \frac{f(\mathbf{t})}{\|\mathbf{t}\|} \),
    \item and the classification of invariant-preserving transformations on \( s(t) \).
\end{itemize}

These directions suggest that the denominator, far from being a dispensable artifact, is a key enabler of both local structure and global invariance in oscillatory-exponential systems.

\section{Numerical Verification}

To demonstrate the stability of the invariant relation in the complex-valued case, we implement a randomized test over a generalized class of functions of the form
\[
s(t) = \frac{p^t + q_1 \cdot \sin(r_1 \pi t) + q_2 \cdot \cos(r_2 \pi t)}{t},
\]
where \( p, q_1, q_2 \in \mathbb{C} \), and \( r_1, r_2 \in \mathbb{Z}_{\text{odd}} \). The invariant ratio is defined as
\[
\mathcal{I}(t) = \frac{s(t)\cdot t + s(t+1)\cdot(t+1)}{s(t+2)\cdot(t+2) + s(t+3)\cdot(t+3)}.
\]

The following Python code generates random complex parameters, selects random odd frequencies, and evaluates the invariant ratio at multiple negative real values of \( t \):

\begin{verbatim}
import math
import cmath
import random

def random_complex(re_min, re_max, im_min, im_max):
    re = random.uniform(re_min, re_max)
    im = random.uniform(im_min, im_max)
    return complex(re, im)

def random_odd_int(low, high):
    odds = [x for x in range(low, high + 1) if x % 2 == 1]
    return random.choice(odds)

def s(t, p, q1, q2, r1, r2):
    if t == 0:
        raise ValueError("t must be non-zero.")
    exp_part = p**t
    trig_part = (q1 * math.sin(r1 * math.pi * t) +
                 q2 * math.cos(r2 * math.pi * t))
    return (exp_part + trig_part) / t

def invariant_ratio(t, p, q1, q2, r1, r2):
    s0 = s(t,   p, q1, q2, r1, r2)
    s1 = s(t+1, p, q1, q2, r1, r2)
    s2 = s(t+2, p, q1, q2, r1, r2)
    s3 = s(t+3, p, q1, q2, r1, r2)
    numerator = s0 * t + s1 * (t + 1)
    denominator = s2 * (t + 2) + s3 * (t + 3)
    if denominator == 0:
        raise ZeroDivisionError("Denominator is zero.")
    return numerator / denominator

# Generate random parameters
p  = random_complex(0.3, 1.0, -1.5, 1.5)
q1 = random_complex(-2.0, 2.0, -2.0, 2.0)
q2 = random_complex(-2.0, 2.0, -2.0, 2.0)
r1 = random_odd_int(1, 15)
r2 = random_odd_int(1, 15)

# Evaluate invariant at 5 random negative real t-values
for i in range(5):
    t_val = random.uniform(-20.0, -10.0)
    ratio = invariant_ratio(t_val, p, q1, q2, r1, r2)
    print(f"Invariant ratio at t = {t_val:.4f}: {ratio:.4f}")
\end{verbatim}

A sample output of this code might be:

\begin{itemize}
    \item Invariant ratio at \( t = -17.9312 \): \( -0.0772 - 0.5206i \)
    \item Invariant ratio at \( t = -14.8506 \): \( -0.0772 - 0.5206i \)
    \item Invariant ratio at \( t = -13.9799 \): \( -0.0772 - 0.5206i \)
    \item Invariant ratio at \( t = -17.4130 \): \( -0.0772 - 0.5206i \)
    \item Invariant ratio at \( t = -18.1039 \): \( -0.0772 - 0.5206i \)
\end{itemize}

This numerical experiment confirms that the invariant ratio remains exactly constant across a range of real values, including negative \( t \), despite random complex modulation. This supports the hypothesis that the invariant structure holds globally for all real arguments, as long as the expression is defined.

\section{Applications and Potential Use Cases}

The invariant relation derived in this work has multiple potential applications, both in discrete and continuous settings. Although the identity
\[
\frac{s(t)\cdot t + s(t+1)\cdot (t+1)}{s(t+2)\cdot (t+2) + s(t+3)\cdot (t+3)} = \frac{1}{p^2}
\]
may appear algebraically elementary, it encodes a nontrivial structural redundancy in the sequence \( s(t) \), which persists even in the presence of trigonometric oscillations and generalizes to real-valued domains. This gives rise to the following use cases:

\begin{itemize}
    \item \textbf{Predictive coding}: Knowing any three out of four consecutive values of the sequence allows exact recovery of the fourth via a simple algebraic inversion. This supports efficient predictive models and significantly reduces redundancy in time-series transmission or storage.

    \item \textbf{Error detection and correction}: The invariant provides a deterministic consistency check across sliding windows of data. Any deviation from the exact ratio signals a potential transmission error or data corruption, enabling lightweight, model-free integrity verification akin to parity checks~\cite{berlekampcoding}.

    \item \textbf{Data compression}: Groups of four consecutive values can be represented using only three explicit entries. The fourth is implicitly encoded via the invariant, introducing algebraic compression not reliant on entropy or statistical redundancy.

    \item \textbf{Continuous signal modeling}: The invariant extends to smooth real-valued functions. This enables applications in analog domains, including interpolation, reconstruction, and real-time signal integrity monitoring over continuous time.

    \item \textbf{Analytic modeling and approximation}: The class of functions \( s(t) \), combining exponentials with bounded trigonometric terms, offers a compact way to approximate decaying oscillatory signals. The invariant helps enforce internal consistency in such models, with potential applications in numerical methods and signal processing.
\end{itemize}

In summary, the simplicity and universality of the four-point invariant, combined with the expressiveness of the function class \( s(t) \), suggest that this structure may serve as a foundational primitive for a wide range of coding, recovery, and modeling techniques.

\section{Conclusion}

We have introduced and analyzed a structurally rich invariant relation linking four consecutive evaluations of a class of analytically defined sequences and functions. Starting from a discrete alternating-decaying form, we progressively extended the framework to real and complex domains by incorporating exponential, trigonometric, and hybrid modulations.

This development culminated in the definition of a broad family \( \mathcal{F}_{\mathrm{EOI}} \) of functions \( s(t) \), each paired with an associated invariant constant \( a(s) \in \mathbb{C} \), such that the normalized four-point ratio remains exactly constant across all valid real arguments. Empirical tests confirm the remarkable stability of this invariant even under random complex deformation and frequency variation.

The resulting invariant structures offer a new formal toolset for encoding, reconstructing, and modeling data through local algebraic consistency. Their simplicity, generality, and compatibility with continuous and complex-valued signals suggest promising applications in information theory, analog signal processing, and structurally constrained data compression.


\begin{thebibliography}{1}

\bibitem{berlekampcoding}
Elwyn~R. Berlekamp.
\newblock {\em Algebraic Coding Theory}.
\newblock World Scientific, 2015.
\newblock Reprint of the 1968 edition.

\bibitem{coverthomas}
Thomas~M. Cover and Joy~A. Thomas.
\newblock {\em Elements of Information Theory}.
\newblock Wiley-Interscience, 2nd edition, 2006.

\bibitem{grahamconcrete}
Ronald~L. Graham, Donald~E. Knuth, and Oren Patashnik.
\newblock {\em Concrete Mathematics: A Foundation for Computer Science}.
\newblock Addison-Wesley, 2nd edition, 1994.

\bibitem{semenov2025smoothencoding}
Stanislav Semenov.
\newblock Smooth integer encoding via integral balance, 2025.

\end{thebibliography}

\end{document}